\documentclass[a4paper,11pt]{article}
\usepackage{authblk}

\usepackage{stix}
\usepackage{hyperref}
\usepackage{amsmath,amsthm}
\usepackage[capitalise, noabbrev]{cleveref}
\usepackage[utf8]{inputenc}
\usepackage[T1]{fontenc}
\usepackage[ruled]{algorithm2e}
\usepackage{booktabs}
\usepackage{graphicx,graphics,float}

\usepackage{parskip} 
\usepackage{float} 
\usepackage{microtype} 
\usepackage{multirow}
\usepackage{multicol} 

\usepackage{todonotes}

\usepackage[
backend=biber,
isbn=false,
url=false,
maxbibnames=99,
sorting=none,
style=numeric-comp,
]{biblatex}

\renewbibmacro{in:}{}

\renewbibmacro*{issue+date}{%
	\setunit{\adddot\space}%
	\iffieldundef{issue}
	{\usebibmacro{date}}
	{\printfield{issue}%
		\setunit*{\addspace}%
		\usebibmacro{date}}%
	\newunit}

\AtEveryBibitem{%
	\clearfield{pages}%
	\clearfield{note}%
	\clearfield{day}%
	\clearfield{month}%
	\clearfield{translator}%
	\clearfield{urldate}%
	\clearfield{pagetotal}%
	\clearfield{publisher}%
}

\theoremstyle{definition}

\usepackage[section]{placeins}
\makeatletter
\AtBeginDocument{%
   \expandafter\renewcommand\expandafter\subsection\expandafter
    {\expandafter\@fb@secFB\subsection}%
}
\makeatother
\newcommand\blfootnote[1]{%
  \begingroup
  \renewcommand\thefootnote{}\footnote{#1}%
  \addtocounter{footnote}{-1}%
  \endgroup
}

\begin{document}	
\title{Spatio-temporal smoothing and dynamics of different electricity flexibility options for highly renewable energy systems - case study for Norway}
\author[1]{Aleksander Grochowicz\thanks{Corresponding author, \url{aleksgro@math.uio.no}}}
\author[2]{Marianne Zeyringer}
\author[1]{Fred Espen Benth}

\renewcommand{\Affilfont}{\itshape\small}
\affil[1]{Department of Mathematics, University of Oslo, P.O. Box 1053 Blindern, 0316 Oslo, Norway}
\affil[2]{Department of Technology Systems, University of Oslo, P.O. Box 70, 2027 Kjeller, Norway}

\date{\today}
\maketitle

\begin{abstract}
In this article, we investigate the mismatch of renewable electricity production to demand and how flexibility options enabling spatial and temporal smoothing can reduce risks of variability. 
As a case study we pick a simplified (partial) 2-region representation of the Norwegian electricity system and focus on wind power.
We represent regional electricity production and demand through two stochastic processes:
the wind capacity factors are modelled as a two-dimensional Ornstein-Uhlenbeck process and electricity demand consists of realistic base load and temperature-induced load coming from a deseasonalised autoregressive process.
We validate these processes, that we have trained on historical data, through Monte Carlo simulations allowing us to generate many statistically representative weather years.
For the investigated realisations (weather years) we study deviations of production from demand under different wind capacities, and introduce different scenarios where flexibility options like storage and transmission are available.
Our analysis shows that simulated loss values are reduced significantly by cooperation between regions and either mode of flexibility.
Combining storage and transmission leads to even more synergies and helps to stabilise production levels and thus reduces likelihoods of inadequacy of renewable power systems.\blfootnote{
{Accepted version, submitted to Applied Energy. \copyright{} 2023 %
        \hspace{2ex} \href{https://creativecommons.org/licenses/by-nc-nd/4.0}{\centering \includegraphics{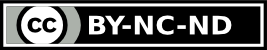}}}%
        \hfill\today
}
\end{abstract}

\newpage

\vspace*{\fill}
\begin{center}
	\textbf{Visual abstract}

	\vspace{2ex}

	\includegraphics{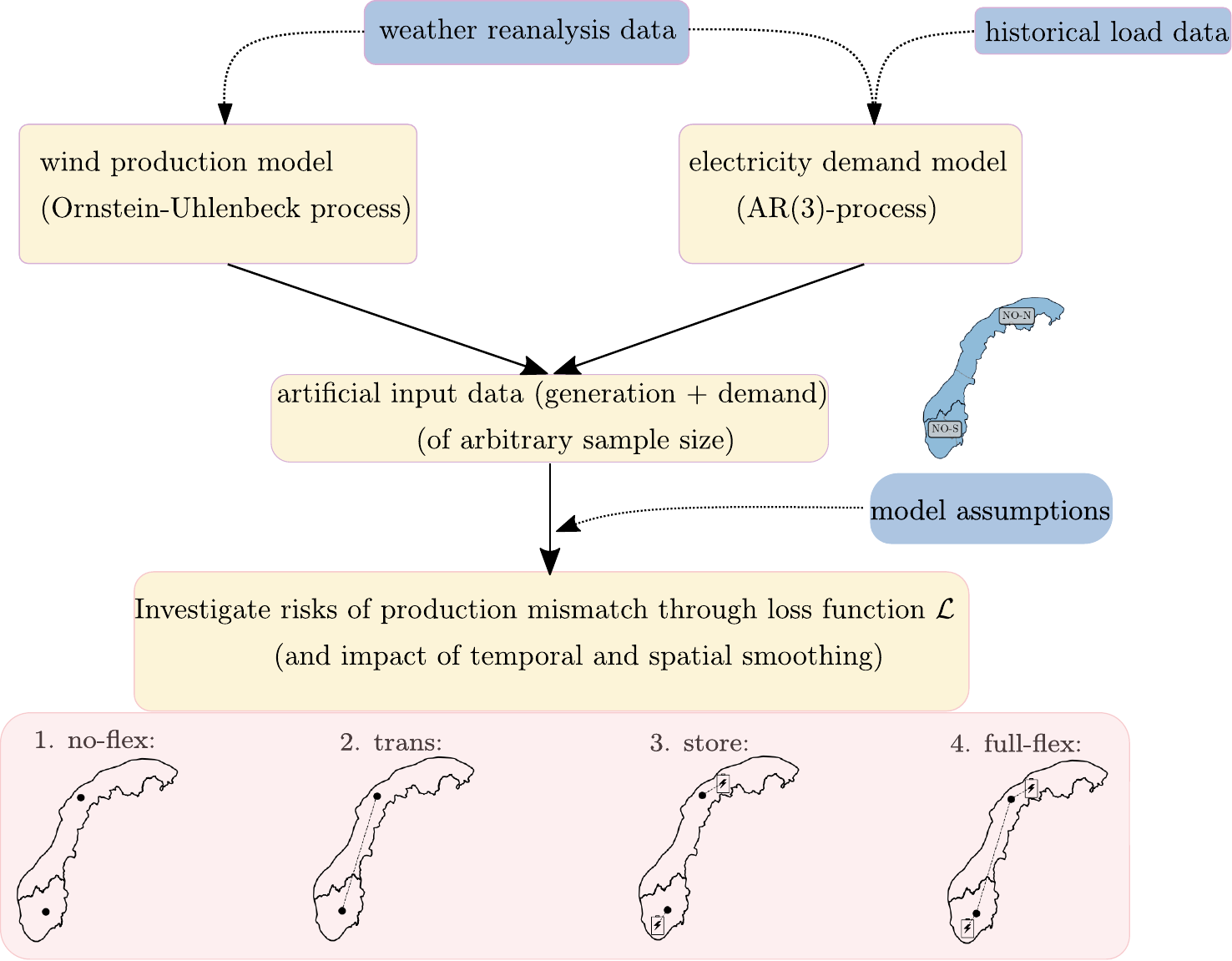}
\end{center}
\vspace*{\fill}

\newpage

\section{Introduction}
\label{sec:intro}
Although virtually all countries signed and ratified the Paris agreement, its implementation is mostly relegated to the national level (even within the European Union). 
In many places the signatories set ambitious targets of decarbonising the power sector \cite{ember-2022} with an increasing penetration of renewable energy sources. 
However, unabated intermittent production can lead to massive mismatch with demand. 
As energy security is a priority of national policy, these matters need to be evaluated through the lens of renewable expansion.
The variability and unpredictability of weather impacts are seen as obstacles, and bracing the power system against undersupply includes incorporating the inherent uncertainties probabilistically. 
However, this is currently out of reach for the multitude of variables and their interactions due to the sheer complexity of weather phenomena. 

In effect, planning and sizing renewable capacities while reducing exposure to weather risks is no easy task and an area of active research, as it is challenging to efficiently include long time series into energy system models \cite{pfenninger-2017}. 
Inter-annual variability is a significant driver of uncertainty for power systems \cite{staffell-pfenninger-2018,zeyringer-price-ea-2018,ringkjob-haugan-ea-2020, grochowicz-vangreevenbroek-ea-2023}, in particular in sizing storage capacities \cite{ruhnau-qvist-2022}.
With increasing renewable penetration, plans for baseload and peak capacities become much more sensitive \cite{bloomfield-brayshaw-ea-2016}.
Since both over- and underproduction are costly and/or inefficient, it is desirable to reduce fluctuations while covering electricity load.
In this context, Olauson et al. study how net load variability in the Nordics can be reduced through different intermittent renewable generation mixes \cite{olauson-ayob-ea-2016}.
The mismatch of demand and renewable generation on a European scale can be captured well by principal components \cite{schwenk-nebbe-vind-ea-2022}, especially for solar \cite{raunbak-zeyer-ea-2017}; in turn, Mühlemann et al. study how meteorological data can be used to site PV production in order to minimise such mismatch \cite{muhlemann-folini-ea-2022}.
Using a Markowitz approach, Roques et al. minimise production variance of wind farm investments \cite{roques-hiroux-ea-2010}.
To compensate the insufficiency of renewable generators, Schlachtberger et al. minimise mismatch by considering necessary back-up dispatch of varying degree of flexibility \cite{schlachtberger-becker-ea-2016}.

Here, we are interested in how flexibility options like transmission and storage contribute to matching demand with intermittent, renewable production, and how their potentials relate to installed capacities. 
Electricity storage allows transferring renewable production in time; transmission allows taking advantage of the different renewable energy potentials across space and that contemporaneous weather conditions can differ from one location to the next. 
We want to evaluate the impact of spatial and temporal smoothing which are necessary to balance higher shares of renewable generation via energy transmission \cite{rodriguez-becker-ea-2014, schlachtberger-brown-ea-2017} and storage \cite{rasmussen-andresen-ea-2012}.
In an IEEE 24-bus power system, Peker et al. study their interaction and synergies through stochastic programming \cite{peker-kocaman-ea-2018}.

There is an increasing number of studies investigating the benefits of storage and transmission to integrate variable renewable energy sources which find that interconnection and storage both decrease costs in a highly renewable US \cite{frew-jacobson-2016} or European power system \cite{bussar-stocker-ea-2016,schlachtberger-brown-ea-2017,child-kemfert-ea-2019,conlon-waite-ea-2019}.
In a more specific context, Bustos et al. illustrate the complementarity of energy storage and transmission expansion for Chile through a mixed-integer linear program \cite{bustos-sauma-ea-2018}. 
Zeyringer et al. show for Great Britain by using 10 weather years that reinforcement of the transmission system consistently lowers system costs \cite{zeyringer-price-ea-2018}. 
Similarly, Roth and Schill find that availability of transmission in a renewable electricity system of Central Europe helps storage to balance wind power and reduces the required capacities of energy storage \cite{roth-schill-2023}. 
Using 41 weather years, Grochowicz et al. demonstrate for a European renewable based power system that transmission expansion is an important contributor to robustness \cite{grochowicz-vangreevenbroek-ea-2023}. 

However, these studies assume perfect foresight and thus do not take into account the stochasticity of the weather variables.  
Here, we add to closing this gap: 
our goal is to investigate the stochasticity of varying weather variables and how storage and transmission levels react to it and quantify the risks of different capacity expansions, based on the flexibility options available. 
Instead of only looking at production variance, we extend this to a system perspective under the inclusion of storage and transmission with the dynamics they introduce, making it necessary to simulate these processes.

Here we answer how spatial and temporal smoothing (or a combination thereof) can mitigate the probability and risk of fluctuations of renewable generation.
We also determine whether introducing these flexibility options influences the sizing potential of wind power capacities and whether benefits of transmission and storage are equally distributed.
These questions can be decisive for public acceptance and adequacy of a transition towards highly renewable power systems.

To accurately quantify the risks of electricity production and demand mismatch, a good representation of production distribution is particularly important for intermittent renewable generation.
In line with this, Sun et al. apply copula theory in order to determine marginal distributions of wind farms in Texas which can enhance forecasting of wind power production \cite{sun-feng-ea-2019}.
In this article we use methods that can generate distributionally consistent artificial data that satisfactorily fit historical data of production and demand based on \cite{benth-pircalabu-2018, benth-christensen-ea-2021}. 
By using artificial data on a daily basis (trained by reanalysis data), we can go beyond historical observations and represent the distribution of weather variables in more detail and through larger samples. 
We then follow a statistical approach to investigate how different flexibility options perform under varying generation capacities that we model in a discretised Ornstein-Uhlenbeck process. 
The performance is measured by a quadratic penalty function on mismatch between electricity supply and demand, signifying that not only are all deviations from target production ``costly'', but larger deviations much more so than slight mismatches: 
overproduction through misinvestment, underproduction through lost load (resulting in blackouts or economic losses through load shedding).
This is demonstrated through a toy-model example with unoperated transmission and storage, inspired by a realistic scenario using current capacity data for Norway.

Norway represents a good case study as it has ambitious goals of electrifying broad swaths of energy sectors, including transport and industry to achieve decarbonisation goals \cite{climate-plan-2022}. 
With this expected increase of electricity load, Norway might transition from being a net exporter to becoming a net importer by the mid-2020s \cite{statnett-kma-2022}. 
Currently there is a lot of dispatch capacity through (reservoir) hydro power (ca. 90\% of production), but through increased integration of variable renewables more and other types of flexibility might be needed, as hydro potential is limited and environmentally contentious \cite{henriksen-2020}.
Additionally, there are disparities in Norway with its five market zones such that the north and south have limited transmission connections and face regional development decisions.
Expansion of wind power --- although not uncontroversial either \cite{linnerud-dugstad-ea-2022,vagero-brate-ea-2023} --- is expected to remain the main intermittent renewable supply source \cite{oed-vindkraft-2022, oed-havvind-2022} as solar resources cannot serve year-round because of the high latitudes.

This article is structured in the following way: 
\cref{sec:methodology} introduces the methodology, background, and implementation while more theoretical details are relegated to the mathematical appendix \ref{sec:math-appendix}.
Afterwards we present our main findings in \cref{sec:results}. Finally we discuss the findings in \cref{sec:discussion} and wrap up with \cref{sec:conclusion}. 

\section{Methodology}
\label{sec:methodology}
\begin{figure}[tbh]
    \centering
    \includegraphics{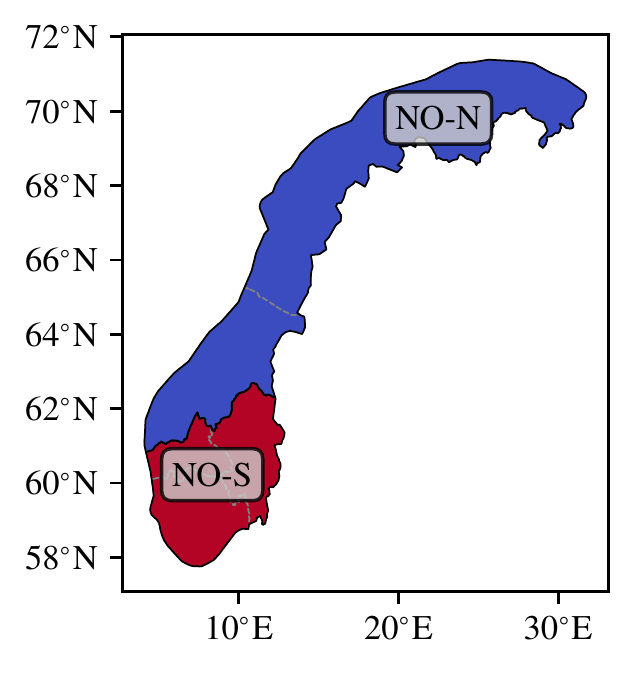}
    \caption{Representation of Norway (which is divided into five bidding zones, NO1-NO5, on the power market) through two regions: ``NO-N'' for northern Norway (NO3 (Trondheim) and NO4 (Tromsø)) and ``NO-S'' for southern Norway (NO1 (Oslo), NO2 (Kristiansand), and NO5 (Bergen)). Currently (2022), 3.25~GW of wind power are installed in NO-N, and 1.8~GW in NO-S.}
    \label{fig:map}
\end{figure}

In this article, we investigate the value flexibility options like transmission and storage represent for renewable power systems through the lens of a quadratic penalty function for mismatch between production and demand (\cref{sec:penalisation}).
Leaning on existing methods, we use two autoregressive stochastic processes to generate artificial input data to obtain larger sample sizes (\cref{sec:artificial-data}).
We do this by representing Norway as a simplified electricity system with two nodes (see \cref{fig:map}); we describe the assumptions in \cref{sec:case-study}.

\begin{figure}
    \centering
    \includegraphics{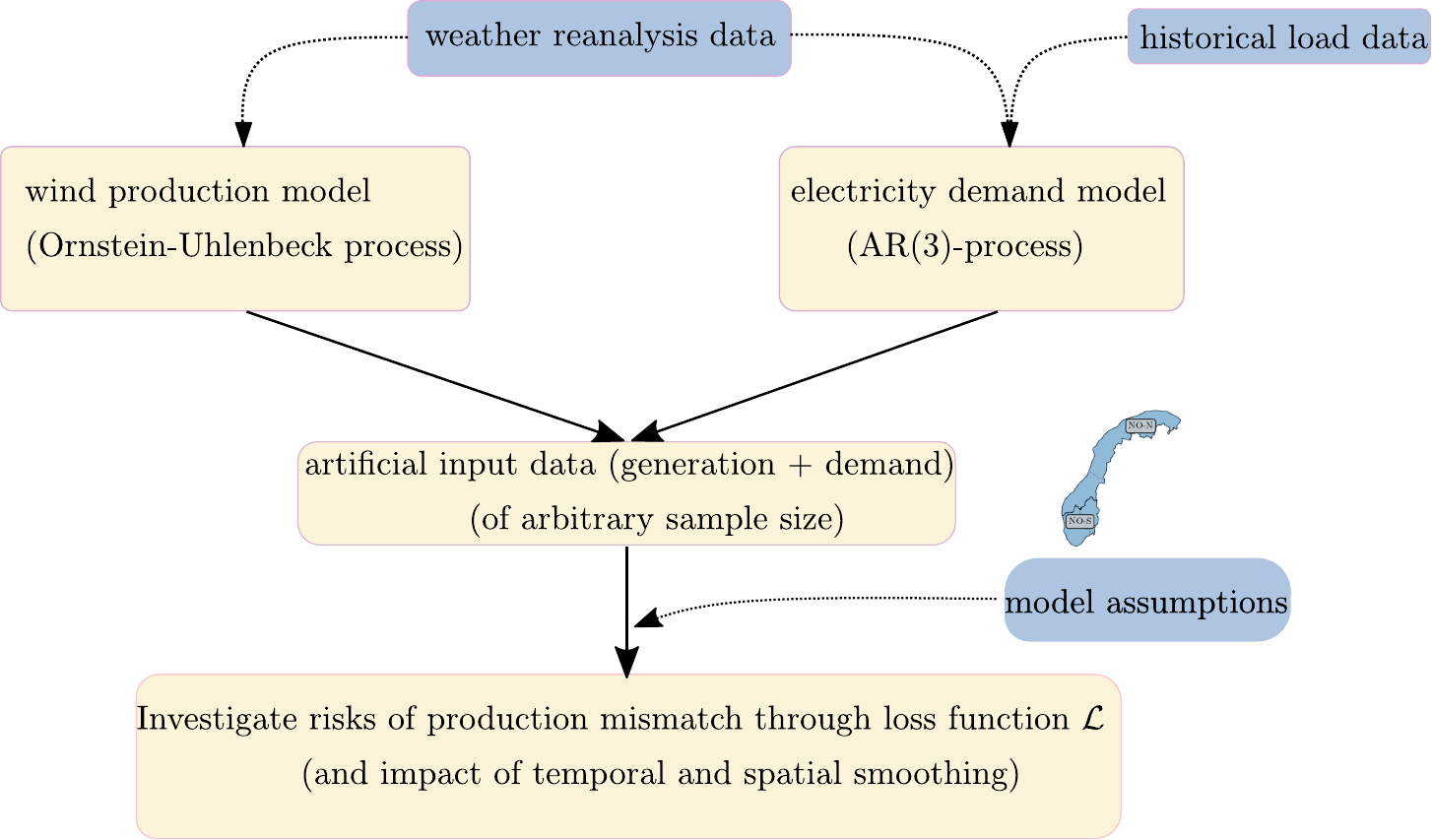}
    \caption{Flowchart describing the procedures and methodology of this article.}
    \label{fig:flowchart}
\end{figure}

\subsection{Penalisation of mismatch between generation and demand}
\label{sec:penalisation}
\begin{figure}
    \centering
    \includegraphics{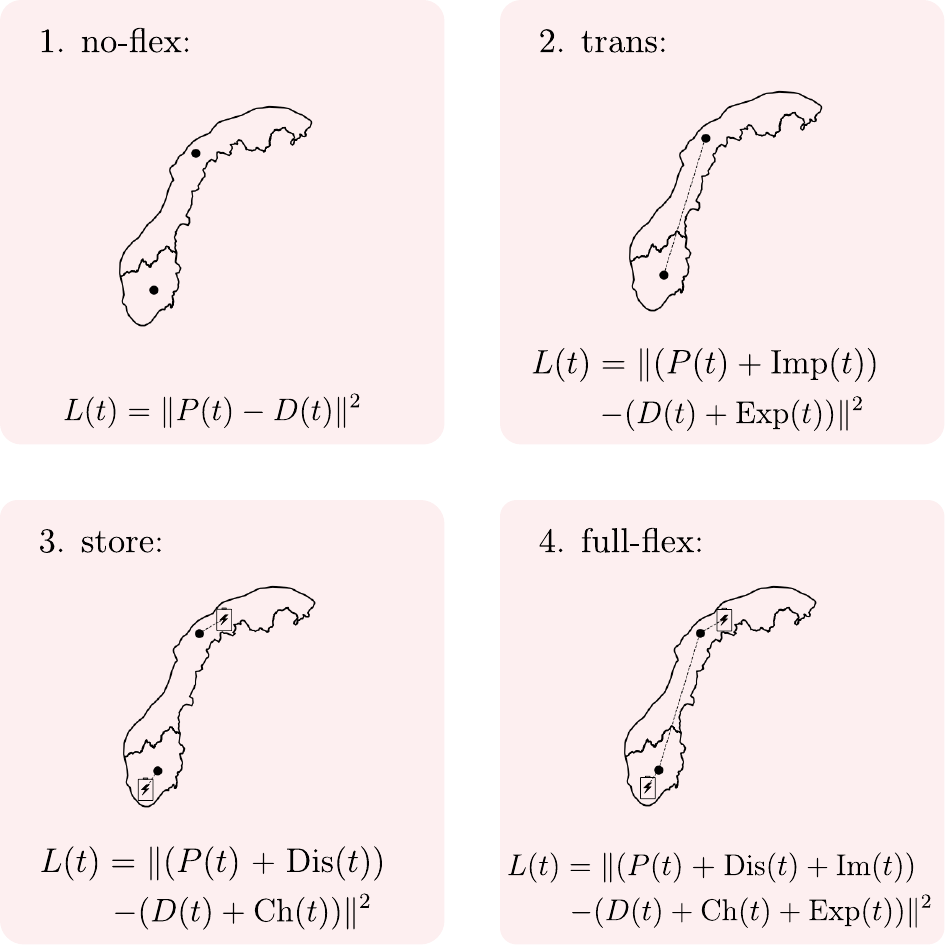}
    \caption{The four different scenarios where we stepwise introduce flexibility through spatial or temporal smoothing. For each scenario, we add the loss functions \cref{eq:loss-simple,eq:loss-trans,eq:loss-stor,eq:loss-flex}.}
    \label{fig:scenarios}
\end{figure}

In this study, we choose wind power production as the generator, thus viewed as a proxy for a renewable system, that needs to cover its share in the electricity production (in daily resolution) consistently throughout the year.
For this we compare four different scenarios (see \cref{fig:scenarios}): one base scenario without flexibility (\emph{no-flex}), then one each with transmission (\emph{trans}) and storage (\emph{stor}) only, respectively, and lastly one with both (\emph{full-flex}).

More precisely, we investigate a quadratic loss function $L(t)$ at time $t$,
\begin{align}
\label{eq:loss-simple}
    L(t) &= \Vert P(t) - D(t) \Vert ^2 \text{ with} \\
    P(t) &= x^T \cdot  C(t), 
\end{align}
where $x \in \mathbb{R}_+^2$ are the installed wind capacities, $C(t) \in [0,1]^2$ the stochastic process describing the regional capacity factors at time $t$\footnote{Here we consider all vectors to be column vectors and $x^T$ is the transpose of $x$.}, and $D(t) \in \mathbb{R}_+^2$ the stochastic nodal demand process at time $t$\footnote{For a more detailed description of the stochastic processes $C$ and $D$ we refer to the mathematical appendix. The capacity factor process here, $C(t)$ is a discretisation of the continuous-time stochastic process $C_t$ we define in \cref{sec:math-appendix}.}. 
Being quadratic, $L$ penalises deviations in both directions, i.e. incorporates trade-offs between under- and overproduction, and does so more the greater the mismatch becomes.
This non-monetary loss function captures the assumption that any deviation from the desired production, i.e. the load at the given time step, is undesirable and that even more so if the deviations increase.

In the scenario \emph{trans}, we connect the two nodes through transmission which allows spatial smoothing by transmitting electricity without efficiency losses to the other node, such that \cref{eq:loss-simple} becomes
\begin{equation}
\label{eq:loss-trans}
    L(t) = \Vert (P(t) + \text{Imp}(t)) - (D(t) +\text{Exp}(t)) \Vert ^2, 
\end{equation}
where Imp$(t)$ and Exp$(t)$ are the imports and exports at time $t$.
For the scenario \emph{stor}, we introduce the possibility to store electricity and thus temporal smoothing; 
to keep capacities comparable, we consider a hypothetical storage (empty at time $0$) unit at each location with the same (dis-)charging capacity as the transmission line.
The dynamics of the storage units $B(t)$ can be described as follows:
\begin{equation}
    \label{eq:dynamics-stor}
    \begin{split}
    \text{Ch}(t) &= \min\{[P(t)-D(t)]^+,B_M-B(t-1), B_C\} \\
    \text{Dis}(t) &= \min\{[D(t) - P(t)]^+,B(t-1), B_D\} \\
    B(t) &= B(t-1) + \text{Ch}(t) + \text{Dis}(t) 
    \text{ and thus } \\
    B(t) &= \sum_{s=1}^t\min\{[P(s)-D(s)]^+,B_M-B(s-1), B_C\} \\
    &-\sum_{s=1}^t\min\{[D(s)-P(s)]^+,B(s-1), B_D \},
    \end{split}
\end{equation}
where $B_M$ are the storage capacities, $B_C, B_D$ the (dis-)charge capacities and $P$ and $D$ are the stochastic production and demand processes as before.

The penalty values for the loss function now become
\begin{equation}
\label{eq:loss-stor}
    L(t) = \Vert (P(t) + \text{Dis}(t)) - (D(t) + \text{Ch}(t)) \Vert^2, 
\end{equation}
where Dis$(t)$ and Ch$(t)$ are the discharge and charge of storage at time $t$.

Finally, in the \emph{full-flex} scenario we can utilise both transmission and storage.
There, overproduction is exported to the node with underproduction (if possible), afterwards locally stored or if possible used to charge the other node's storage (similar to the algorithm in \cite{kaldellis-kapsali-ea-2010}).
Similarly if capacities allow, it is possible for one node with underproduction to discharge and transmit electricity from the other node's storage.
Thus the loss function becomes
\begin{equation}
\label{eq:loss-flex}
\begin{split}
    L(t) &= \Vert (P(t) + \text{Dis}(t) + \text{Im}(t)) 
    - (D(t) + \text{Ch}(t) + \text{Exp}(t)) \Vert^2.
\end{split}
\end{equation}
Note that in all these cases, the flexibility options are reactive (we have no foresight) and cannot be operated and must follow the described rules.

The loss functions \cref{eq:loss-simple,eq:loss-trans,eq:loss-stor,eq:loss-flex} are related to the objective function used in \cite{schlachtberger-becker-ea-2016} where the authors combine a quadratic penalisation of underproduction with minimising installed capacities and maximising usage of back-up capacities.
By penalising overproduction as well, we implicitly discourage overcapacities although we systematically study the interactions of flexibility options and their usage under different capacities.

\subsection{Artificial input data to represent stochasticity of weather}
\label{sec:artificial-data}
To investigate the stochasticity of varying weather variables and how storage and transmission levels react to it, we quantify the risks of different capacity expansions, based on the flexibility options available.
For this, we generate daily capacity factors for wind power as introduced by Benth et al. \cite{benth-christensen-ea-2021} to obtain many statistically representative weather years.
These data are synthesised from ERA5 reanalysis data ranging from 1980 to 2020 \cite{era5-data} and possess similar moments and distributional qualities.
Furthermore, we estimate mean electricity load for many realisations of a calendar year, driven by an autoregressive (AR) process which is trained on ERA5 reanalysis temperature data of the current representative period 1991--2020.
These data follow realistic daily load profiles obtained through a regression as in \cite{grochowicz-vangreevenbroek-ea-2023} with stochasticity introduced by the AR process describing temperature throughout a weather year \cite{eggen-dahl-ea-2022,benth-benth-ea-2008, hardle-cabrera-ea-2010,hardle-cabrera-2012}. 

Note that we use daily average load, since we use calendar days as time steps.
The approach to modelling demand and production as stochastic processes is described in more detail in \cref{sec:math-appendix}.
The synthetic data allow us to take many samples and consider the (expected) penalty values 
\begin{equation}
\label{eq:penalty}
    \mathcal{L} = \mathbb{E} \left[\sum_{t=1}^{365} L(t) \right],
\end{equation}
where $L(t)$ is defined as in \cref{eq:loss-simple,eq:loss-trans,eq:loss-stor,eq:loss-flex}; we also study the penalty values as a mean process $\mathcal{L}$, depending on different capacities $x \in \mathbb{R}^2$ as in \cref{eq:loss-simple}.

\subsection{Case study of Norway}
\label{sec:case-study}

We represent the Norwegian electricity system as a simplified two-node system and apply the aforementioned methodology.
Such a division of the country is justified by the limited existing transmission capacities between the two parts of the country; in fact, electricity flows to a much larger extent via the north-south connections in Sweden due to much higher capacities.
The transmission capacities between NO-N and NO-S are estimated by the Norwegian TSO to be between 900~MW and 1400~MW (although these flows are rarely reached) \cite{statnett-transport-2021}, whereas the Swedish transmission corridor has current capacities around 7000~MW \cite{statnett-svenskakraftnat-2022}.
Thus we use a transmission capacity of 900~MW for the scenarios \emph{trans} and \emph{full-flex} and for comparison purposes, set this as a value for (dis-)charging capacities $B_C, B_D$.
Following the analyses in \cite{otero-martius-ea-2022,raynaud-hingray-ea-2018} about lengths of wind droughts in Norway, we select storage capacities of 15~GWh in NO-N and 30~GWh in NO-S, just enough to cover 14 time steps of high demand in each region\footnote{Note that we work with daily averages to simplify the relationship with existing capacities.}.
Since Norway is divided into five market zones, investment decisions directly impact the respective market zone and lead to different prices within the country.
This resembles more the situation across borders in Europe than within many (larger) European countries like France or Germany.

Additionally, the southern market zones are connected more strongly to the European power market being home to the majority of the Norwegian interconnectors and all DC connections (more than 7 out of the approx. 9~GW capacity) \cite{hofstad-askheim-ea-2022}, and their prices are more interwoven with those on the European continent \cite{statnett-2018-prices,statnett-cableprice}.

We consider a renewable power system for which we aim to minimise the risk of deviations of power demand to power production, and compare scenarios where we introduce the options of 1. transmitting electricity between the zones, 2. storing electricity, or 3. both combined.
The current Norwegian electricity system is relying mostly on hydropower (91\% of production in 2021 \cite{ssb-electricity-production-2022}) and with its hydro reservoirs (capacity equivalent to 70\% of annual electricity consumption \cite{ssb-electricty-consumption-2022, nve-magasinstatistikk}) much of the energy can be stored and dispatched.
So far, the installed pumped hydro storage is far below its potential and furthermore not a realistic consideration for this purpose, as it is designed to be operated seasonally during flood season \cite{pitorac-vereide-ea-2020,egging-tomasgard-2018}.
Still, we assume efficiencies of 75\% for charging and 90\% for discharging as in \cite{PyPSAEur} for pumped hydro storage; however, we consider this storage to be generic (so no inflow despite using efficiency assumptions of pumped hydro).

Due to decarbonisation and electrification targets, the country is expected to become a net importer by the mid-2020s after being a significant exporter for most years in the last decades \cite{statnett-kma-2022}.
Thanks to favourable energy resources, wind energy has expanded from less than 1~GW to 5~GW in 5 years and expected to be further expanded both on land and sea \cite{oed-vindkraft-2022, oed-havvind-2022}.
The currently installed wind capacities (2022) are 3257~MW in NO-N and 1811~MW in NO-S and have an annual average generation of 16.8~TWh \cite{nve-windturbines}, constituting 12.8\% of the average electricity demand for the years 2014--2018 in Norway, based on Nordpool data \cite{nordpoolgroup-2021}.
Thus we set the demand process $D(t)$ which we use for this analysis to cover 12.8\% of the estimated Norwegian load\footnote{After adding market zone demand to the regions, the average Norwegian annual electricity demand for NO-N is 43.7~TWh and 87.8~TWh for NO-S \cite{nordpoolgroup-2021}, so we are aiming to cover 5.6~TWh in NO-N and 11.2~TWh in NO-S on average with wind power.} at all times. 
We disregard the remaining bulk of Norwegian load, assuming it is covered by other production means (such as hydro) and that these other technologies cannot cover the demand in our analysis.

\section{Results}
\label{sec:results}
We present the results from the analysis we conduct based on \cref{sec:methodology}. 
First, we state general observations of the loss surfaces we obtain for the four different scenarios.
For each of the scenarios, we look at the capacities for which the expected penalty value is lowest and explore these further to obtain distributional properties of the simulated processes. 

Finally, we compare the characteristics of how storage and transmission influence the mismatch across scenarios.

We compute the loss functions as described in \cref{sec:methodology} for 100 realisations, i.e. 100 distinct weather years, and thus estimate $\mathcal{L}$ as in \cref{eq:penalty} for wind capacity combinations between 3.25~GW (current) and 6~GW in NO-N and 1.8~GW (current) and 10~GW in NO-S.
The current capacities are thus always included in these combinations.
It should be noted that for these ranges, there can be both over- and undercapacities, which both lead to higher penalty values of $\mathcal{L}$.
We observe that in all simulations, any mode of flexibility always improves the penalty value compared to the scenario without flexibility.
However, the effectiveness of the flexibility modes varies and depends on the installed capacity of wind power generation, as shown by the optimal generation capacities (in the sense of minimising $\mathcal{L}$) in \cref{tab:opt_caps}.
As capacity factors tend to be higher in NO-N (see \cref{fig:cap_factors}), the availability of transmission shifts production away from NO-S (see \emph{full-flex} vs. \emph{stor}) such that electricity flows primarily from north to south.
\begin{table}[htb]
    \centering
    \begin{tabular}{c|c|c|c|c}
        Scenario & Description & NO-N & NO-S & Improvement \\ \toprule
        Current & Minimal capacities for all scenarios. & 3.25~GW & 1.8~GW & 0\%  \\ \midrule
        no-flex & Generation expandable. & 3.25~GW & 4.7~GW & 36\%  \\ \midrule
        trans & Transmission available; generation expandable. & 3.25~GW & 3.8~GW & 55\% \\ \midrule
        stor & Local storage available; generation expandable. & 3.25~GW & 6.45~GW & 65\% \\ \midrule
        full-flex & Scenarios stor and trans combined. & 3.25~GW & 5.55~GW & 84\% \\ \bottomrule
    \end{tabular}
    \caption{Optimal wind capacities and improvement of loss values (by scenario compared to current capacities (2022)). All scenarios only consider expanding generation capacities. Storage and transmission capacities are fixed as in \cref{sec:methodology}, a sensitivity analysis is conducted in \cref{sec:sens-analysis}.}
    \label{tab:opt_caps}
\end{table}

\begin{figure}[htb]
    \centering
    \includegraphics{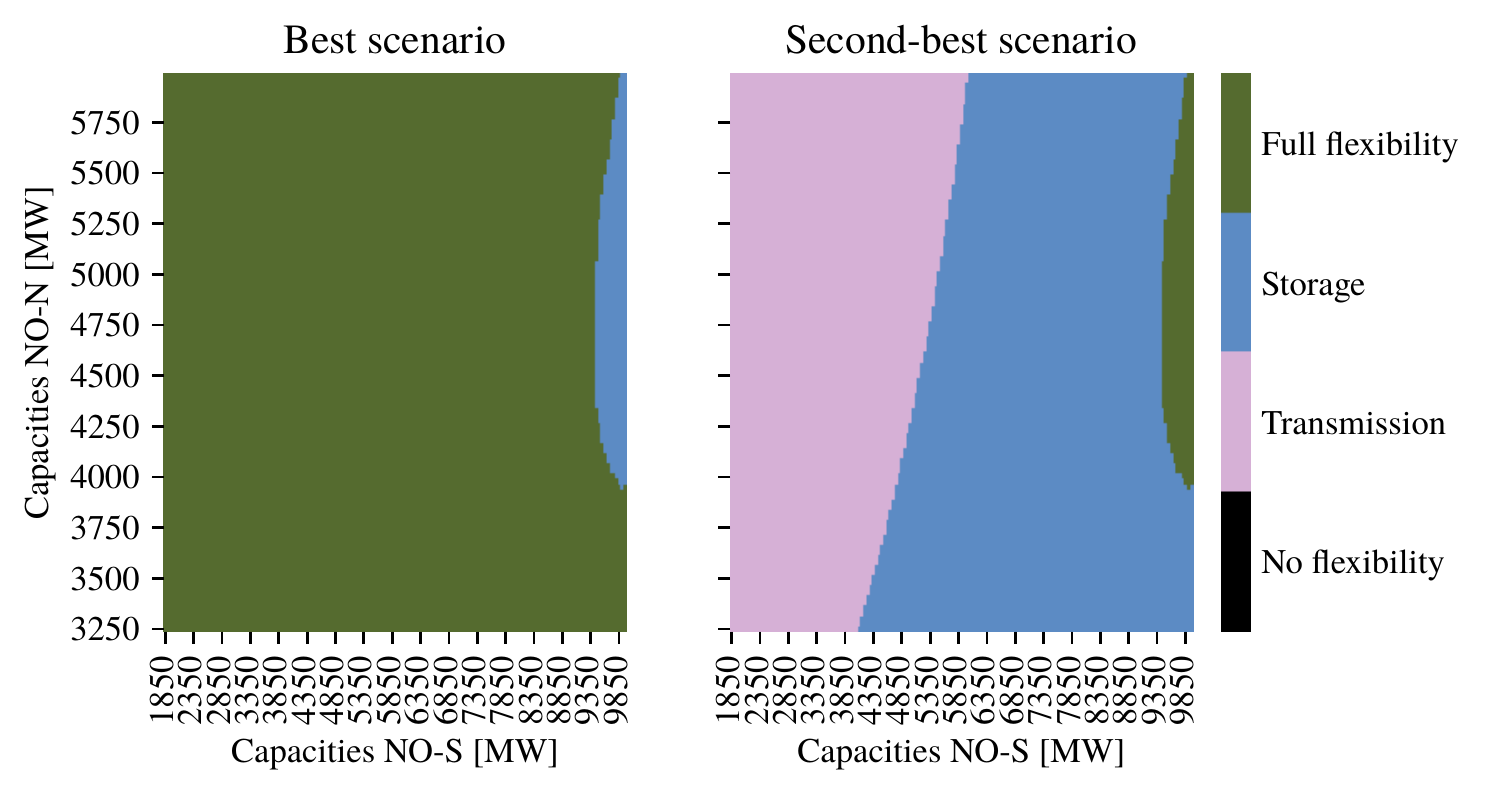}
    \caption{Ordering of the best two scenarios by wind power capacities in NO-N and NO-S. Left: the scenario with the lowest penalty for given capacities is marked; right: the scenario with the second-lowest penalty for given capacity is marked.}
    \label{fig:surfaces}
\end{figure}

Since the load that is supposed to be covered by wind power relates to current production, hefty expansion of generation capacities leads to overcapacities, which is why optimal capacities for \emph{no-flex} do not deviate from the currently installed ones in NO-N (being currently much higher than in NO-S, despite lower load). 
The same values for their respective loss functions could be achieved with lower capacities for the remaining scenarios;  it is possible to increase the capacities beyond these levels to actually reduce the penalty values by up to 84\%. 
\Cref{fig:surfaces} shows that for almost all capacity combinations, having both transmission and storage reduces the penalty values the most; only at the extremes (right corner) the overcapacities are so large that also having transmission \emph{increases} the penalty values (due to higher round-trip efficiency than storing and discharging). 

At the same time, there is a clear curve where the benefits of storage outweigh those of transmission.
This is where the loss surfaces --- which are convex and paraboloid-shaped because of the quadratic loss function --- for the two scenarios intersect: in practice, transmission is more valuable than storage at lower capacities in general, especially when those in NO-N are higher.
Still, the benefit of both transmission and storage does outweigh the sum of the components for a significant part of the studied combinations.

\subsection{Different scenarios}
\subsubsection{Status quo}
\begin{figure}
    \centering
    \includegraphics{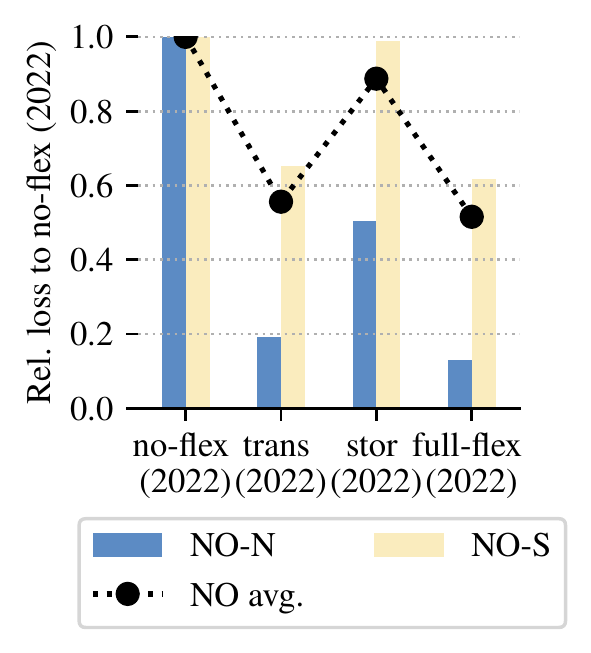}
    \caption{Relative improvement of penalty values with current capacities (2022) through introduction of flexibility, by region.}
    \label{fig:current_scenarios}
\end{figure}

At the current (2022) capacities, the benefits of flexibility in terms of the loss function are limited due to the intermittency of production against the (relatively) consistent profiles of load --- recall that their total sum are equal on average. 
Most of the capacities are installed in northern Norway, where load is lower, which leads to overcapacities there and undercapacities in southern Norway in relation to designated load.
Therefore the scenarios with transmission allow exports from north to south, which reduces penalty values significantly in both regions.
If storage is available, this allows some temporal smoothing, however this is limited by the capacities as shown in \cref{fig:comp_level} (which is equivalent for NO-N); capacities and thus production tend to be too low for any build-up of storage in NO-S. 
Regardless, all flexibility scenarios lead to improvements, ranging from 52\% with storage and 87\% with both flexibility options in northern Norway, and from 1\% with storage to 37\% with both flexibility options in southern Norway (see \cref{fig:current_scenarios}).

\subsubsection{No flexibility}
\begin{figure}[htb]
    \centering
    \includegraphics{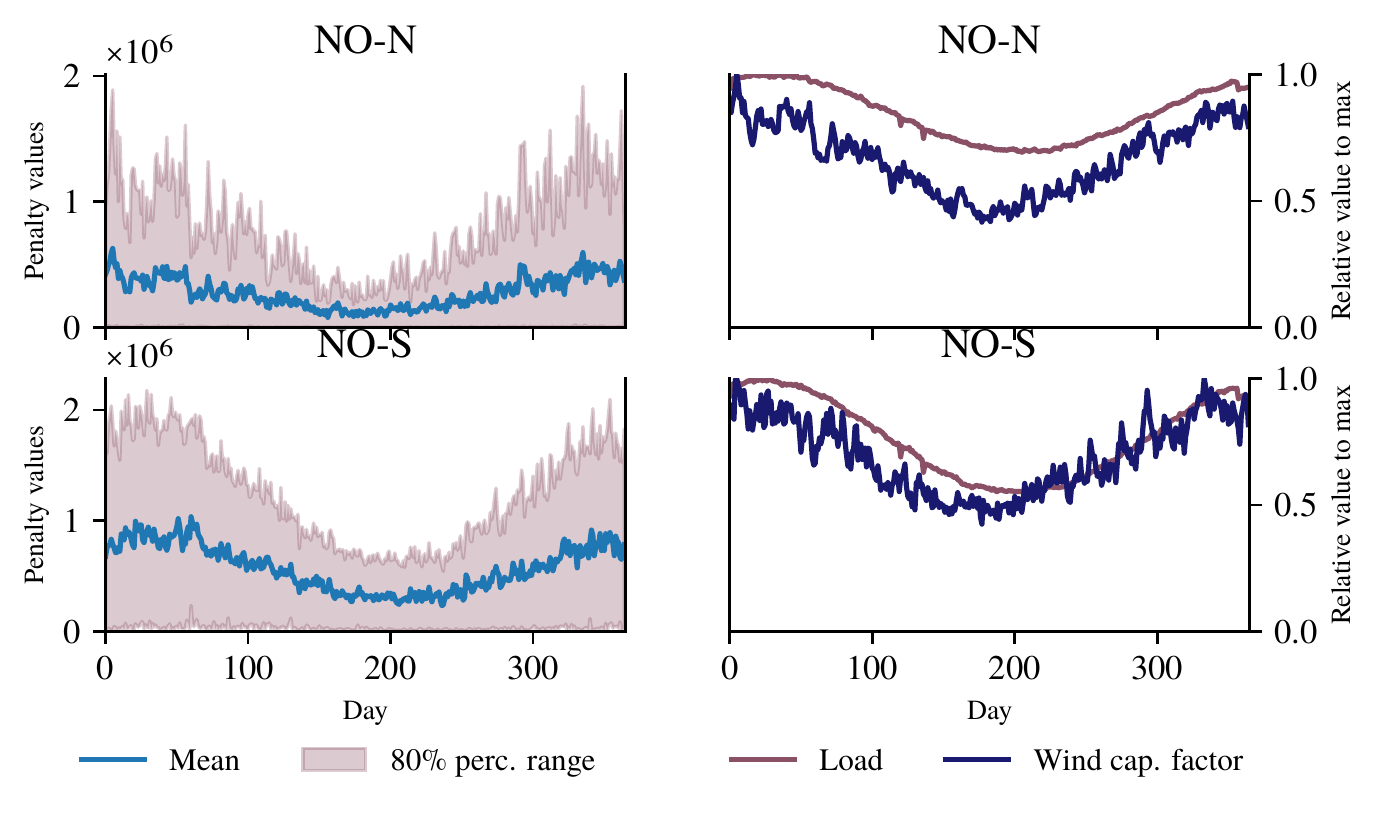}
    \caption{Left: Average penalty values with probability ranges for the scenario without flexibility. Right: Seasonal profiles of average load and average wind capacity factors (normalised). 
    For both, the average values are divided by the corresponding maximum value.}
    \label{fig:no-flex}
\end{figure}

In the \emph{no-flex} scenario without any flexibility, the loss values can be reduced by additional wind capacities in the south of Norway (from 1.8 to about 4.7~GW) as \cref{fig:comp_penalties} depicts, while capacities and therefore penalty values in NO-N remain the same. 
\cref{fig:no-flex} shows that the largest penalties occur in the winter months in both regions, with larger uncertainties and erratic behaviour in NO-N (based on 100 simulations), as both under- and overproduction can lead to mismatch.
In NO-S, the seasonal profiles are much more pronounced: penalty values and their variability decrease sharply.
The smaller bands of probability can be explained by a general decrease of demand in the summer months: possible mismatches for low wind periods are thus smaller.
However, the seasonal profiles of average wind capacity factors and relative load are more similar in southern Norway which explains the reductions in the summer (see \cref{fig:no-flex}).

\subsubsection{Flexibility through transmission}
\begin{figure}[htb]
    \centering
    \includegraphics{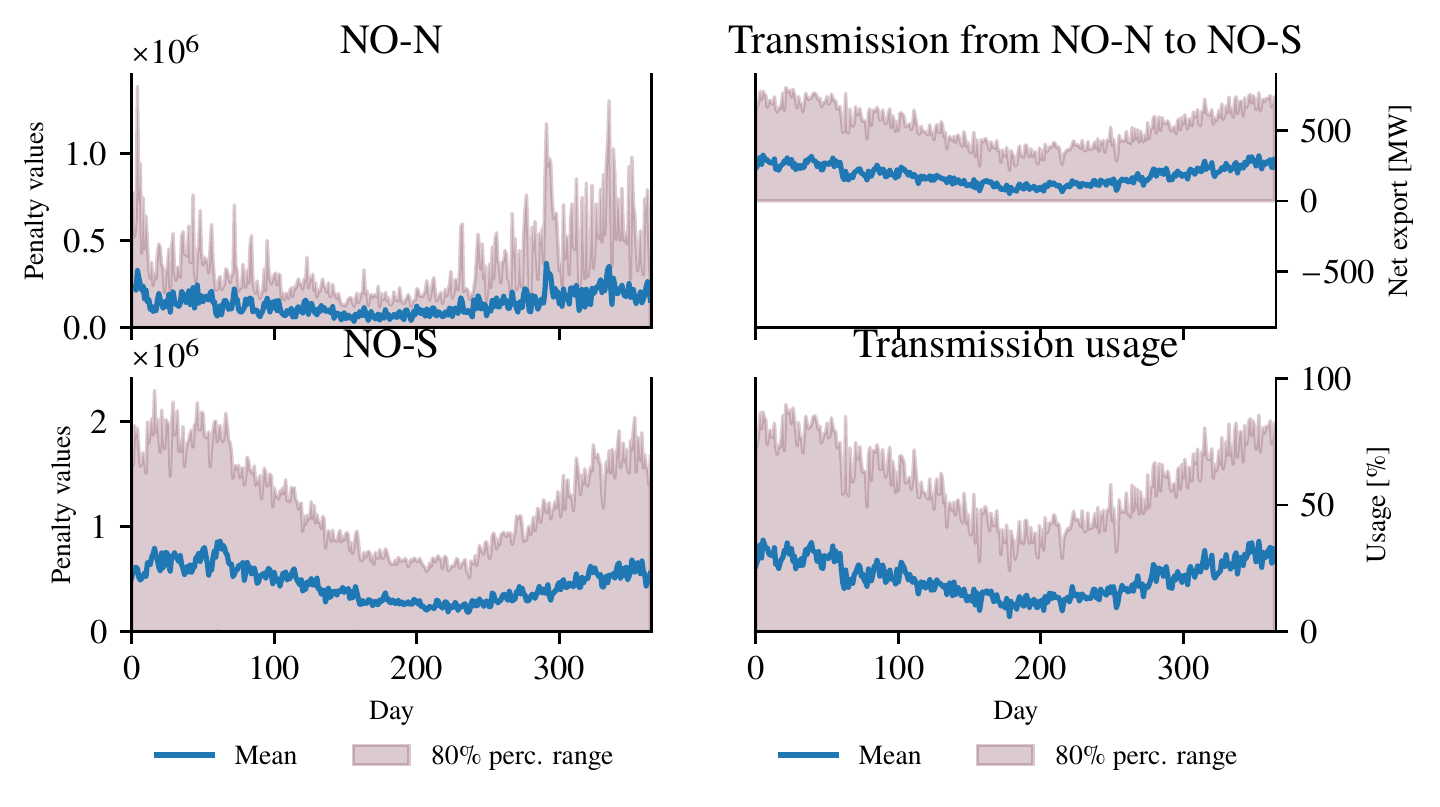}
    \caption{Left: Average penalty values with probability ranges for the scenario with transmission. Right top: average net exports from NO-N to NO-S with probability ranges; right bottom: average transmission usage with probability ranges}
    \label{fig:trans_penalties}
\end{figure}

In the \emph{trans} scenario, we introduce transmission of up to 900~MW between NO-N and NO-S.
With this first flexibility option, the loss-minimising capacities in NO-S are higher than the current ones (as in \cref{tab:opt_caps}) but lower than \emph{no-flex} scenario, in fact exactly by the capacity of transmission, which appears to be mostly unidirectional.
\cref{fig:trans_penalties} shows that the loss function for NO-S has lower values through the introduction of transmission and is flatter than before.
In NO-N, the penalty values (especially the quantiles) are much more volatile than before, but on average reduced meaningfully.
Due to higher capacities in the north of the country, the more extreme values come from realisations with high capacity factors and lead to significant overproduction (which is also penalised by the loss function).
Especially in the winter months, there is more export from NO-N to NO-S (see \cref{fig:trans_penalties}) than in the summer, but the flow is not limited by the capacities.

\subsubsection{Flexibility through storage}
\begin{figure}[htb]
    \centering
    \includegraphics{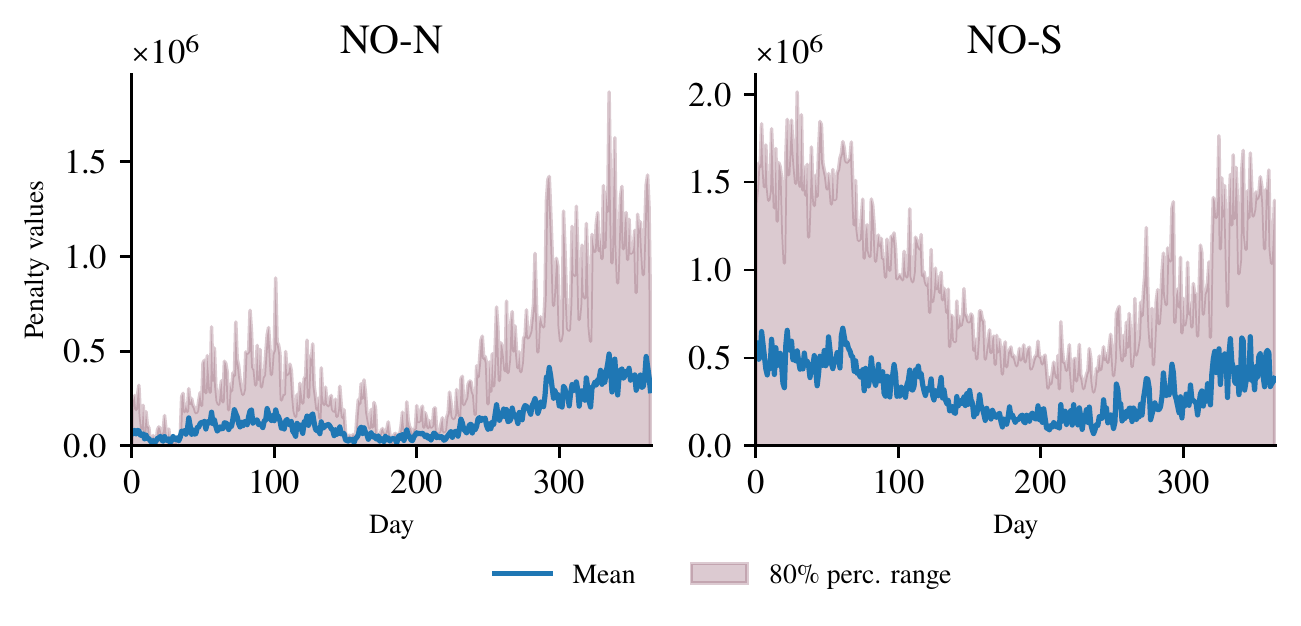}
    \caption{Penalty values with 80\% probability ranges with optimal capacities for the scenario with storage as the only source of flexibility.}
    \label{fig:stor_penalties}
\end{figure}

\begin{figure}[htb]
    \centering
    \includegraphics{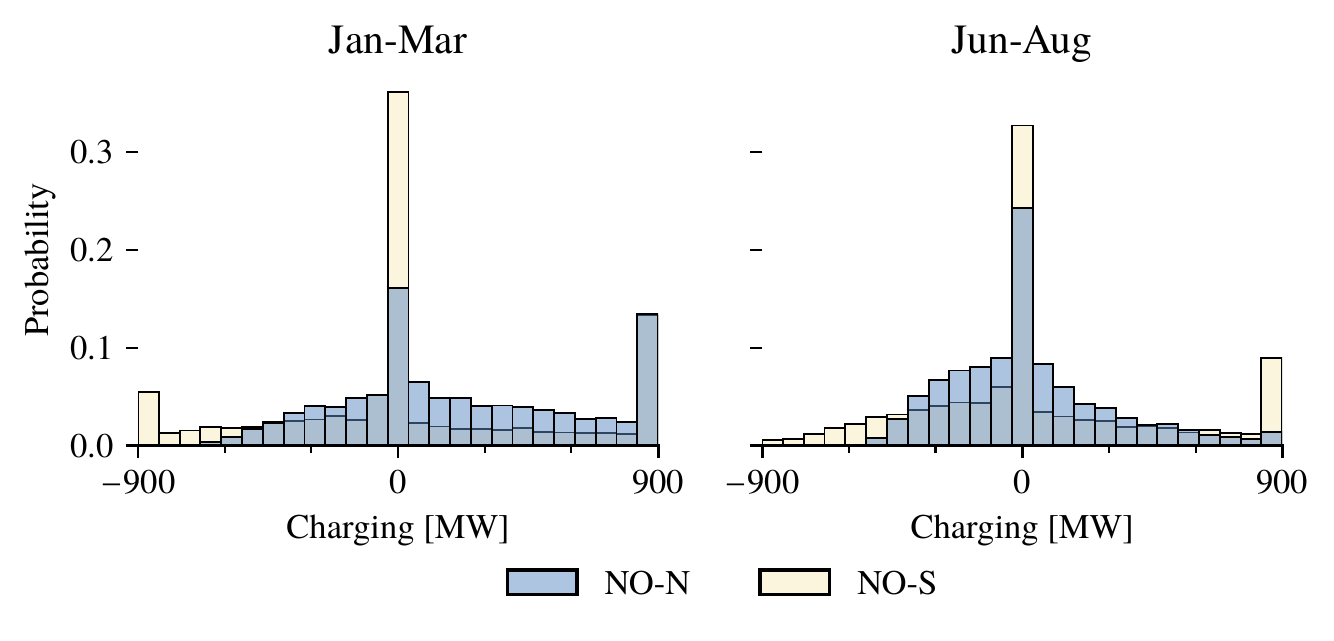}
    \caption{Histogram of storage usage by region in the winter (Jan-Mar) and summer months (Jun-Aug). A high accumulation at maximal charging (900~MW) or discharging capacity (-900~MW) points at insufficient capacities.}
    \label{fig:stor_usage-distr}
\end{figure}

\begin{figure}[tbh]
    \centering
    \includegraphics{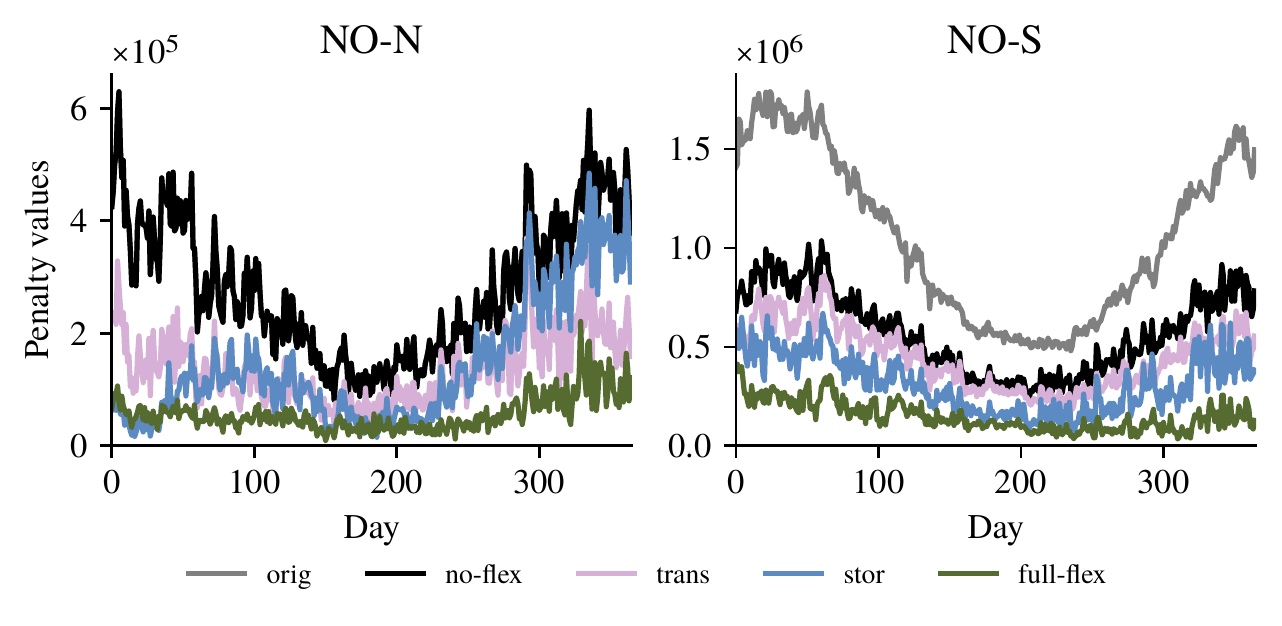}
    \caption{Average penalty values depending on the scenario. Note the different order of magnitude between NO-N and NO-S.}
    \label{fig:comp_penalties}
\end{figure}

The possibility to store overproduction in the \emph{stor} scenario, leads to much larger capacities, of which 6.45~GW are allocated in NO-S and, unchanged, 3.25~GW in NO-N.
The penalty values are similar to \emph{trans}, with some improvements in NO-S due to the optimal capacities and higher penalties in NO-N because of overproduction that exceeds storage capacities in the spring already (as in \cref{fig:comp_penalties} and \cref{fig:stor_penalties}).
These penalties in northern Norway decrease once the wind becomes weaker in the summer months (see \cref{fig:no-flex}), more so than the load profile.
This suggests that already for the existing capacities in NO-N, a larger fraction of electricity load than what we assume could be covered by wind power.
Large variations in southern Norway can still occur, especially early in the year, when storage needs to be built up. 
On the other hand, there appears to be less overproduction in NO-S, where average storage levels are low (as shown in \cref{fig:comp_level}) with a 30-60\% chance of the storage being empty throughout the year (see \cref{fig:stor_level-distr}).
\cref{fig:stor_usage-distr} shows that charging and discharging capacities in NO-S are not large enough, especially in the winter months; these limitations are less of a concern in NO-N, perhaps because storage levels remain close to their capacity limit throughout the year.

\subsubsection{Full flexibility through transmission and storage}
When both transmission and storage are available, the capacities in NO-S (while those in NO-N remain the same throughout the scenarios) are lower than in the \emph{stor} scenario (by 900~MW, which coincides with the transmission capacity), but higher than in the other scenarios.
The loss function is reduced significantly, mostly through a strong decrease in NO-S (as these values were several times higher than those in NO-N). 
This again is possible through overproduction in NO-N which can be stored in many realisations or sent southward.
Due to slightly higher usage of transmission than before (mostly from NO-N to NO-S), these two modes of flexibility work symbiotically.
As previously observed, the charging capacities can be a limitation (see \cref{fig:full-flex_usage-distr}), being more pronounced in the winter.
The storage levels show similar trends as in the \emph{stor} scenario, however thanks to transmission, the storage capacities in northern Norway are not filled up as early; in southern Norway average levels are similar throughout the year (despite lower generation capacities) and with a more uniform distribution.
After January, the likelihood of storage being empty in NO-S lies still between 30-60\% and can reach 35\% in NO-N in July.
Taking both regions into account, the integration is highly beneficial, reducing the value in \cref{eq:loss-flex} by 75\% in comparison to the \emph{no-flex} scenario (47\% to \emph{stor} and 65\% to \emph{trans}).

\begin{figure}[htb]
    \centering
    \includegraphics{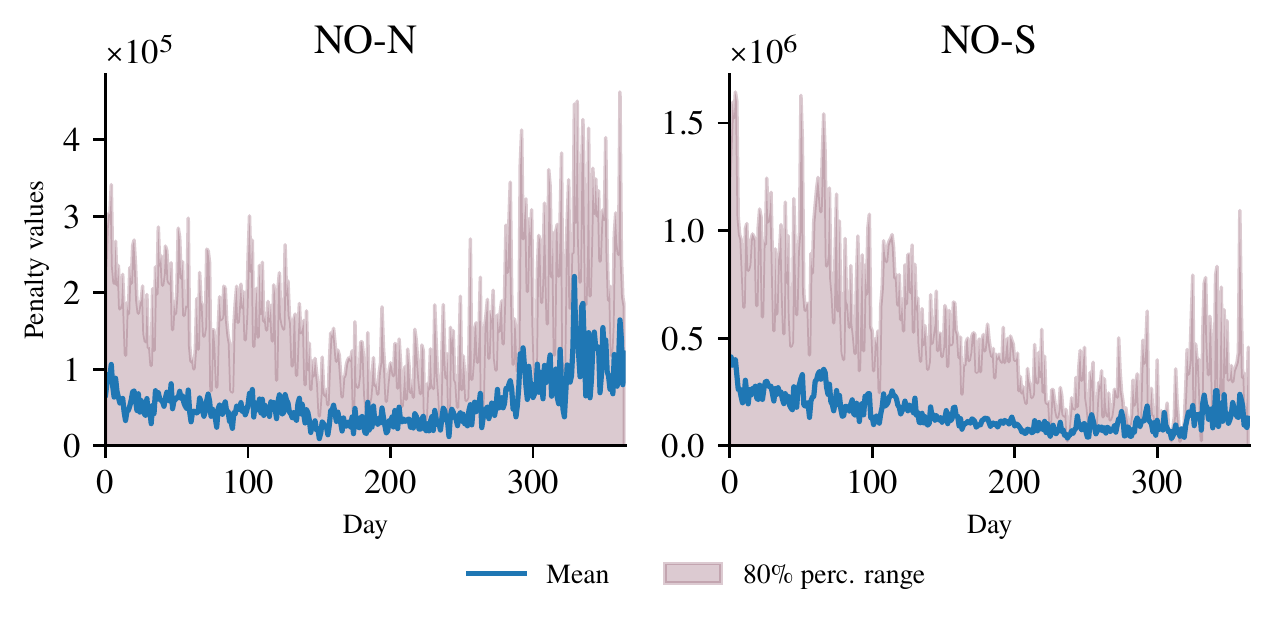}
    \caption{Penalty values with 80\% probability ranges with optimal capacities for the scenario with storage and transmission as sources of flexibility.}
    \label{fig:full-flex_penalties}
\end{figure}

\subsection{Comparison of the scenarios}

Due to the seasonal profiles and variance of wind capacity factors and demand profile in Norway, the penalty values in the summer are lower than in the winter when demand is higher and mismatch can be higher. 
The only exception occurs in the \emph{stor} scenario in northern Norway, where penalty values (and their variability) increase once storage is filled up and overproduction leads to losses.

In general, lower wind speeds in the summer follow the lower demand quite well, although the average capacity factors in NO-N decrease more than load, where heating demand persists also throughout the summer. 
If storage is available, it often comes from overproduction early in the year and is used more in the second half of the year.
Note that minimising the overall penalty value could have adverse effects for one node as the aggregate values are dominated by NO-S, although this does not take the stored electricity into account (e.g. in NO-N) whose value is not reflected in the loss function.
With the availability of transmission (in the \emph{trans} and \emph{full-flex} scenarios) the better wind conditions in northern Norway encourage the surplus of additional capacities to be exported.
The ability to store electricity benefits both nodes, in particular with higher capacities in northern Norway where good wind conditions early in the year produce electricity to be stored into the more difficult summer.
It remains to stress that in this study the integration of as much flexibility as possible and cooperation is beneficial for both net importer and exporter (apart from total welfare benefits).

\section{Discussion}
\label{sec:discussion}
We simulate wind power production and electricity demand as stochastic processes, the former being a discretised version of an autoregressive process, the latter being based on a regression and another autoregressive process.
We use these simulations in a simplified representation of Norway with two zones investigating how flexibility provided by transmission and storage impact the mismatch between production and demand.
In particular we can assess the benefits of spatial smoothing (through transmission) with those of temporal smoothing (through storage), depending on the installed capacities.

We find that any mode of flexibility helps to better match demand and production while also reducing the risk of large deviations --- which here can prove ``costly'' both as overproduction as well as underproduction.
For the purpose of this study, overproduction is penalised equally to underproduction (in the quadratic loss function); however, in highly renewable systems the mismatch does in fact inflict costs through overcapacities, inefficient usage of back-up or baseload capacities or even opportunity costs through curtailment.
The interaction of transmission with storage changes the dynamics and alters some distributional qualities of the stochastic process describing storage.
Low or negative correlation between nodes can make spatial smoothing through transmission even more valuable.
In general, as seasonal patterns of production and demand match, there is a higher risk of mismatches in winter as variability has a larger impact than the higher average production.
As the variability weighs heavier on the penalty function, this impacts optimal capacities more than in usual capacity expansion problems where installation is driven by periods of lowest production.

Even with modest flexibility, we show that in our study both regions benefit from investing in both storage and transmission, but to a different extent due to uneven distribution of capacities.
Concretely, Northern Norway has more favourable wind conditions than the southern part, and thus has export potential that disproportionally benefits Southern Norway.
This raises questions of distributive justice, and is a relevant aspect in domestic politics where national climate targets (such as electrification of gas fields \cite{officeofthenorwegianprimeminister-2023}) and electrification come in direct conflict with local electricity demand and land use of indigenous people \cite{thesupremecourtofnorway-2021}.
The perceived fairness of such a situation may also affect the socio-political feasibility of increased (spatial) smoothing \cite{wolsink-2007,vuichard-broughel-ea-2022}.
Mutual benefits (also for exporters) should thus be a priority of political questions for integration and collaboration to be successful and socio-politically accepted.

Although we consider a reality-inspired application, there are several deviations from a realistic use case in our simplified set-up.
For instance, we set over- and underproduction as equally undesirable, which explains higher penalties in the \emph{stor} scenario for NO-N once storage capacities are filled up.
At the same time, our modelling of storage (which complicates the simulations significantly in contrast to transmission) is purely reactive in order to understand how the stochastic process is affected by the interactions with transmission.
Our implementation without foresight or operational decision stands in contrast to reality where there is some foresight and short-term weather forecasting is improving.
Similarly the scope of this study is limited insofar as its spatial (two nodes for Norway) and temporal (daily) resolution is fairly limited, as we focus more on the availability of many realisations of the processes and a systematic evaluation of flexibility options under different capacities.
Since we choose to simulate a large grid of capacities ($110 \cdot 163$ combinations), we choose 100 realisations/weather years for this study due to high computational burden (the simulations took approx. 3-4 days on a 64 cores machine) --- with more realisations it is likely that results are smoother.
Each realisation is based on a calendar year and could be shifted to fit storage patterns more realistically. 
Still, 100 realisations constitute a larger sample than we can obtain from most time series (such as weather reanalysis data), and the here presented way of generating artificial weather data can be valuable for risk analyses of renewable generation.

\section{Conclusion and outlook}
\label{sec:conclusion}
In this article, we evaluate how transmission and storage (in isolation and a combination thereof) affect stability of renewable generation.
A simplified representation of the Norwegian electricity system in which we describe electricity demand, generation, and storage through discretised continuous-time stochastic processes shows that combining flexibility options and cooperation is beneficial and allows a higher penetration of renewable technologies.
Model formulations that focus on minimising variability can accompany capacity expansion problems in energy system optimisation models that, in contrast to this study, take monetary considerations into account.
Comparing these perspectives could help to understand possible trade-offs of different objectives.

Although we focus in this particular study on Norway, certain insights can be extended beyond Norway.
Without much adaptation, such a study can also be conducted in other parts of the world with one dominant renewable technology which matches seasonal demand well.
Extensions to other technologies or mixes of technologies are also possible, but might require different stochastic models (e.g. for solar irradiation).
The market design in Norway with numerous different zones and transmission bottlenecks makes the situation comparable to other regions in Europe and the US;
as national investment decisions do for Europe, regional decisions can impact adequacy, costs and prices for the different market zones in Norway.
At the same time, these methods, and in particular the synthesis of artificial data, could be valuable to assess and hedge investment risks in renewable generation portfolios.

Lastly, it would be interesting to include more locations and modelling a more realistic transmission network, but also including more variables such as e.g. solar irradiation (and correlate them) and more technologies in order to obtain a more realistic set-up which was deemed out of scope for this case study.

\section*{Acknowledgements}

We would like to thank the anonymous reviewers for their helpful comments and suggestions. 
We thank Oskar Vågerö for his inputs on socio-political acceptance of wind power in Norway.

The computations were partly performed on the Norwegian Research and Education Cloud (NREC), using resources provided by the University of Bergen and the University of Oslo. \url{https://www.nrec.no}

\printbibliography

\appendix
\section{Mathematical appendix}
\label{sec:math-appendix}
\subsection{Load data synthesis}
We generate synthetic daily load data for the two parts of Norway based on historical consumption data which we divide into temperature-dependent and -independent components.
We model temperatures in Norway by a deseasonalised autoregressive (AR) process based on the current reference period (1991 -- 2020) of temperature reanalysis data \cite{era5-data}.
With these two processes we can generate time series of realistic load, based on stochasticity of temperature, which we now summarise briefly.

For this, we use ERA5 reanalysis temperature data \cite{era5-data} (averaged over the grid cells in each region) and load data for the five market zones from Nordpool (for the years 2014 to 2018) \cite{nordpoolgroup-2021} which we aggregate to the two nodes, NO-N (consisting of NO3 and NO4) and NO-S (consisting of NO1, NO2, NO5).

We follow the second part of the regression from \cite{grochowicz-vangreevenbroek-ea-2023} in which we estimate the daily load $D_i(t)$ of different weekdays (and holidays which are classified as Sundays) that is independent of temperature and the load that is mostly driven by heating demand (as cooling demand is statistically insignificant in the north and not that relevant in the south): 
\begin{align}
\begin{split}
     D_i(t) &= \beta_{\text{weekday}, i}(t) + \beta^{\text{heating}}_i \cdot \text{HDD}_i(t) + \beta^{\text{cooling}}_i \cdot \text{CDD}_i(t) \\ 
     \text{  with  } \text{HDD}_i(t)  &= \max\{15.5 - T_i(t), 0\} \text{ and } \text{CDD}_i(t)  = \max\{T_i(t) - 15.5, 0\},
\end{split}
\end{align}
where $\beta_{\text{weekday}, i}$ are the regression parameters for the temperature-independent load for different weekdays and $\beta^{\text{heating}}_i$, $\beta^{\text{cooling}}_i$ are the regression parameters for one degree of heating/cooling demand (measured by HDD/CDD, the heating/cooling-degree day index) for each region $i$.
Having estimated the temperature-independent base load, what is left for a realistic time series is to introduce heating demand which constitutes a major part of electricity demand in Norway and cooling demand.

As we want to simulate many possible demand years, we introduce randomness and variability through an appropriate model of temperature data.
We describe daily average temperatures as a deseasonalised AR(3)-process (similar to \cite{eggen-dahl-ea-2022, benth-benth-ea-2008}), which allows us to generate an arbitrary number of realisations of demand years as follows.
Based on ERA5 reanalysis temperature data aggregated to NO-N and NO-S for 1991 -- 2020, we can describe the temperature process $T(t)$ as a sum of a seasonal $S(t)$ and deseasonalised component $Y(t)$: 
\begin{equation}
    T(t) = S(t) + Y(t) = (a + b \cdot \sin(\frac{2 \pi t}{365}) + c \cdot \cos(\frac{2 \pi t}{365})) + Y(t),
\end{equation}
where we estimate the parameters $a, b, c \in \mathbb{R}$.
Afterwards we estimate the partial autocorrelation function for $Y(t)$, noticing that $Y(t)$ now captures the temperature anomalies at time $t$ from the seasonal mean at that time step. 
From the partial autocorrelation function we infer a lag of $3$ which explains the choice of the AR(3)-model\footnote{In many of the above-cited papers on temperature models, continuous-time autoregressive models are considered, see e.g \cite{benth-benth-ea-2008,hardle-cabrera-ea-2010,hardle-cabrera-2012}.} whose parameters we also estimate.
From the seasonal component $S(t)$ and simulating the AR(3)-process we now generate 100 weather years of temperature data --- more precisely, the heating degree days (HDD) --- which we use together with the temperature-independent base load and temperature-dependent parameters from the above regression to obtain realistic load time series for NO-N and NO-S.

\subsection{Capacity factor synthesis}

We build on the approach by Benth et al. in \cite{benth-christensen-ea-2021} for generating daily wind capacity factors for 100 weather years in our case. 
They find that it is possible to represent wind capacity factors for different locations as a multi-dimensional Ornstein-Uhlenbeck process with L\'evy jumps.
For clarity, we summarise and extract the formulations we use in this paper in the following.

Suppose that $d$ is the number of locations (here $d=2$) and $X_t$ is a $d$-dimensional Ornstein-Uhlenbeck (OU) process of the form
\begin{equation}
\label{eq:ou-process}
dX_t=-\Lambda X_t dt+\Sigma dL_t
\end{equation}
where $L$ is a $d$-dimensional compound Poisson process with independent exponential jumps.
We assume $\Lambda =\text{diag}(\lambda_1, \dots, \lambda_d)$ to be a diagonal matrix describing the lag and $\Sigma=(\sigma_{ij})_{i,j} \in \mathbb{R}^{d \times d}$ describing the idiosyncratic risks from $L$.
The values $\lambda_i$ > 0 and $\sigma_{ij}$ > 0 are strictly positive.

This OU process is connected to the stochastic process representing the capacity factors $C_t \in \mathbb{R}^d$ as follows:
\begin{equation}
    C_t = (1 \dots 1)^T - \exp(-S_t \cdot X_t),
\end{equation}
where 
$$
S_t=\left[\begin{array}{cccc} s_{1}(t) & 0 & \cdots & 0   \\
0 & s_{2}(t) & \cdots & 0  \\
\cdot & \cdot & \cdots & \cdot  \\
0 & 0 & \cdots & s_d(t) 
\end{array}\right]\in\mathbb R^{d \times d}
$$ 
consists of seasonality functions $s_i(t) = a_i + b_i \cdot \sin(\frac{2 \pi t}{365}) + c_i \cdot \cos(\frac{2 \pi t}{365})$ for each location $i$.

\begin{figure}[htb]
    \centering
    \includegraphics{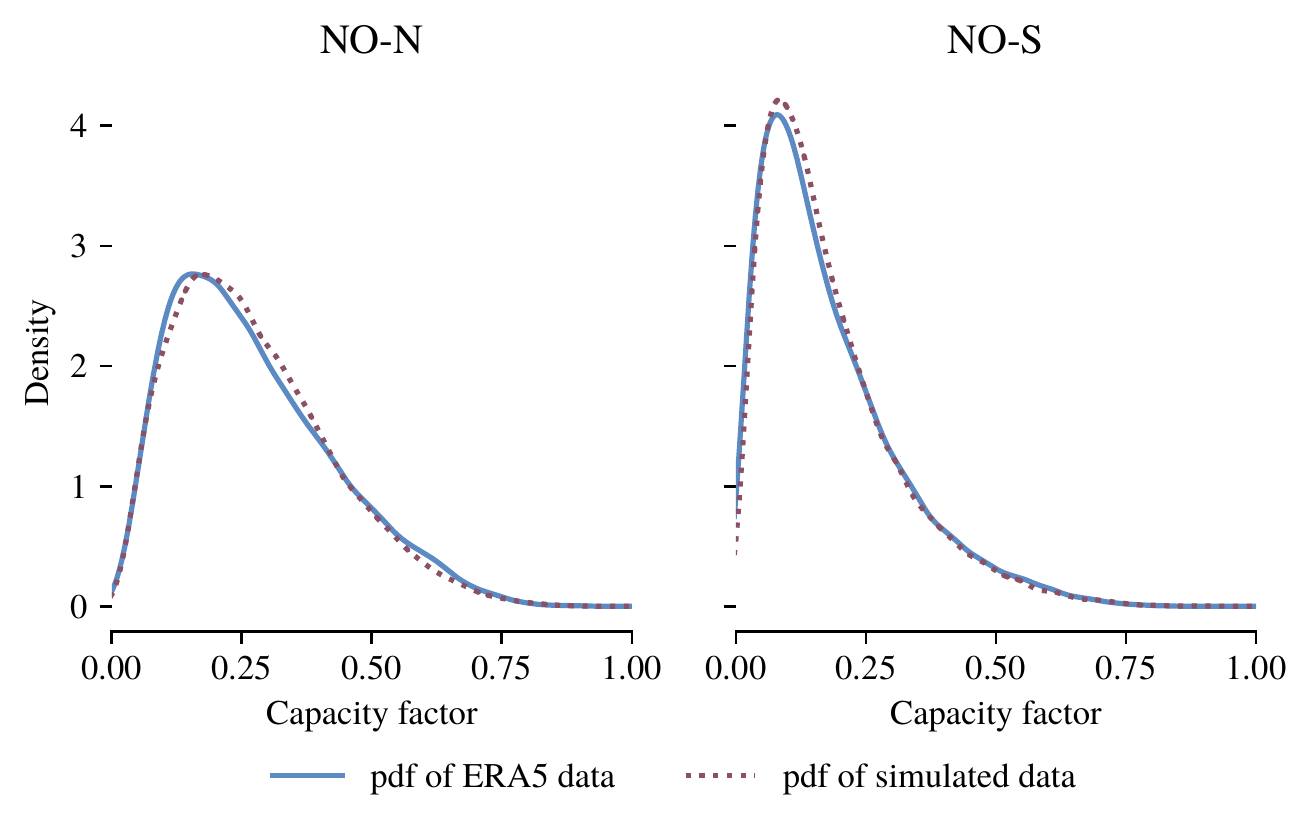}
    \caption{Probability distribution functions of the original ERA5 capacity factors (solid) for northern and southern Norway and the capacity factors $C$ obtained through the Ornstein-Uhlenbeck process $X$ (dotted).}
    \label{fig:cap_factors}
\end{figure}

\begin{figure}[htb]
    \centering
    \includegraphics{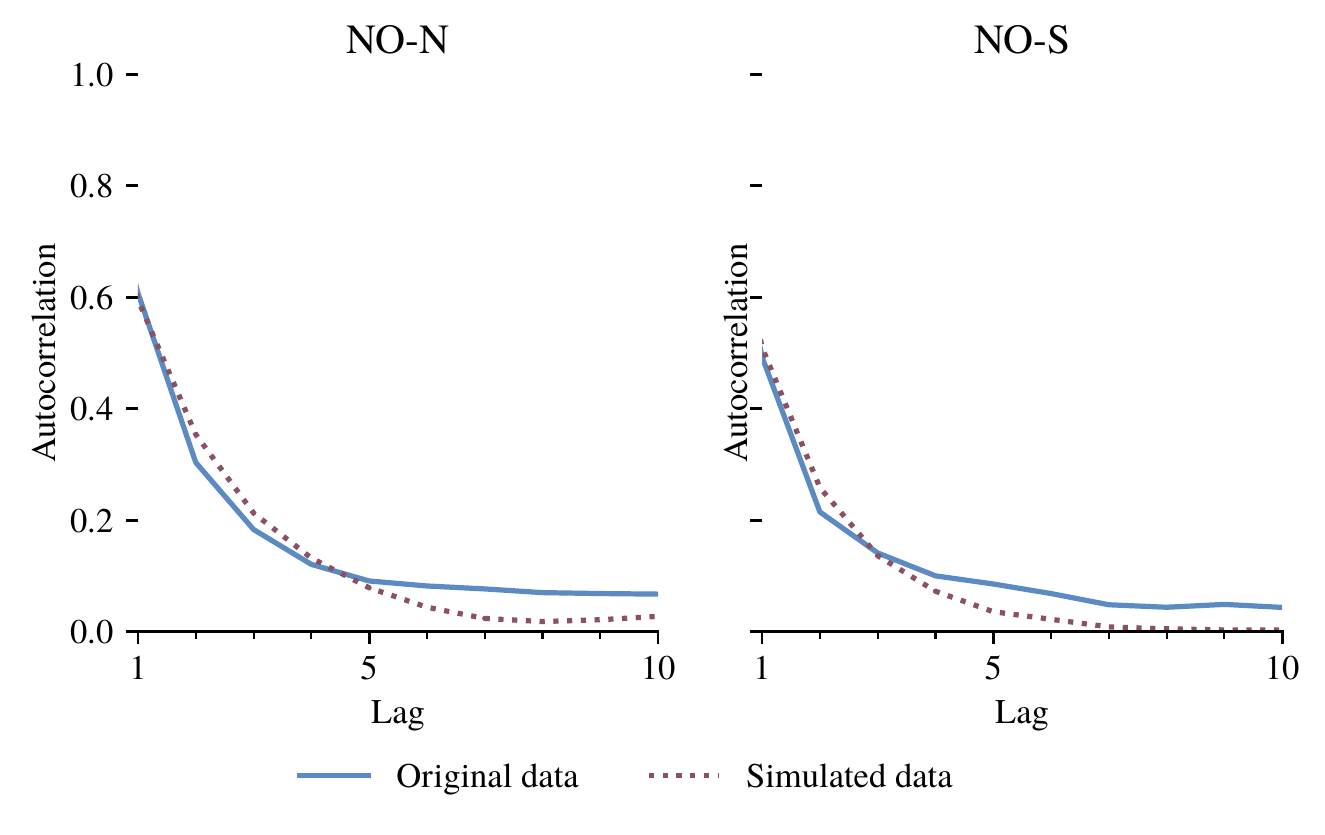}
    \caption{Autocorrelation of the original process $X$ estimated from ERA5 data and the simulated process.}
    \label{fig:autocorrelation}
\end{figure}

We discretise the continuous-time stochastic process, and consider each day of the year as a time step. This we do in order to estimate the model to the following data. 
We use reanalysis data from ERA5 \cite{era5-data} from 1980 until 2020 and the open source tool Atlite \cite{hofmann-hampp-ea-2021} to compute capacity factors for northern Norway (defined by the market zones NO3 and NO4) and southern Norway (NO1, NO2, and NO5).
With these capacity factors, we follow the steps described in \cite{benth-christensen-ea-2021} to estimate  $\Lambda, \Sigma, S_t$ and $L_t$.
Using these estimates, we obtain a stochastic time-series model fitted to data whose simulation matches the first three moments and autocorrelation of the original capacity factors as shown in \cref{fig:cap_factors,fig:autocorrelation} and \cref{tab:stat-caps}.
A finer and more sophisticated autoregressive model than \cref{eq:ou-process} could further reduce the mismatch for higher lags and improve model quality.
At the same time, our simulations show a good fit in the correlation (ERA5: 0.464, OU: 0.470) between capacity factors in NO-N and NO-S.

\begin{table}[htb]
    \centering
    \begin{tabular}{c|cc|cc}
         & ERA5: NO-N & OU: NO-N & ERA5: NO-S & OU: NO-S \\ \toprule
        Mean & 0.273  & 0.269 & 0.180 & 0.180  \\ \midrule
        Standard deviation & 0.155 & 0.149 & 0.135 & 0.131 \\ \midrule
        Skewness & 0.752 &  0.767 & 1.251 & 1.362 \\ \bottomrule
    \end{tabular}
    \caption{Mean, standard deviation, and skewness for the original capacity factors (ERA5, left) and the synthetic capacity factors generated by the OU process (OU, right).}
    \label{tab:stat-caps}
\end{table}

We generate $100$ calendar years of wind capacity time series on a daily basis, but with this approach, we can also generate much larger numbers, probably smoother results, and rely more on sampling than we do here (in particular to obtain the surfaces in \cref{fig:surfaces} based on the average of the same 100 years).  

\section{Sensitivity analysis}
\label{sec:sens-analysis}
We conduct a sensitivity analysis on four components, varying
\begin{itemize}
    \item demand ($\pm 10\%$ jointly in both regions, $\pm 10\%$ for each single region)
    \item transmission capacity ($\pm 10\%$, $\pm 50\%$)
    \item storage capacity ($\pm 10\%$ jointly, $\pm 50\% $ jointly)
    \item (dis-)charging capacity ($\pm 10\%$ jointly, $\pm 50\%$ jointly)
\end{itemize}
and studying how that impacts optimal capacities and the loss functions.

The impact of transmission capacities on generation capacities and penalties is rather small, only when transmission is halved (reduction by 450~MW) do we observe increased installations in NO-S (by 200~MW for \emph{trans} and 400~MW for \emph{full-flex}).

Variation of storage capacities influence penalty values --- but not capacities --- for \emph{stor}: as storage is filled up by spring in NO-N, increasing storage capacities helps to delay overproduction reducing loss values.
In the situation where both storage and transmission are available, increased storage leads to a shift of capacities from south to north (by 200~MW with 150\% capacity) such that more overproduction can be stored and exported. 
This lowers penalty values in NO-S (and keeps them constant in NO-N) such that the actual benefits are reaped in the importing region and not where additional capacity is installed.

\begin{figure}[tbh]
    \centering
    \includegraphics{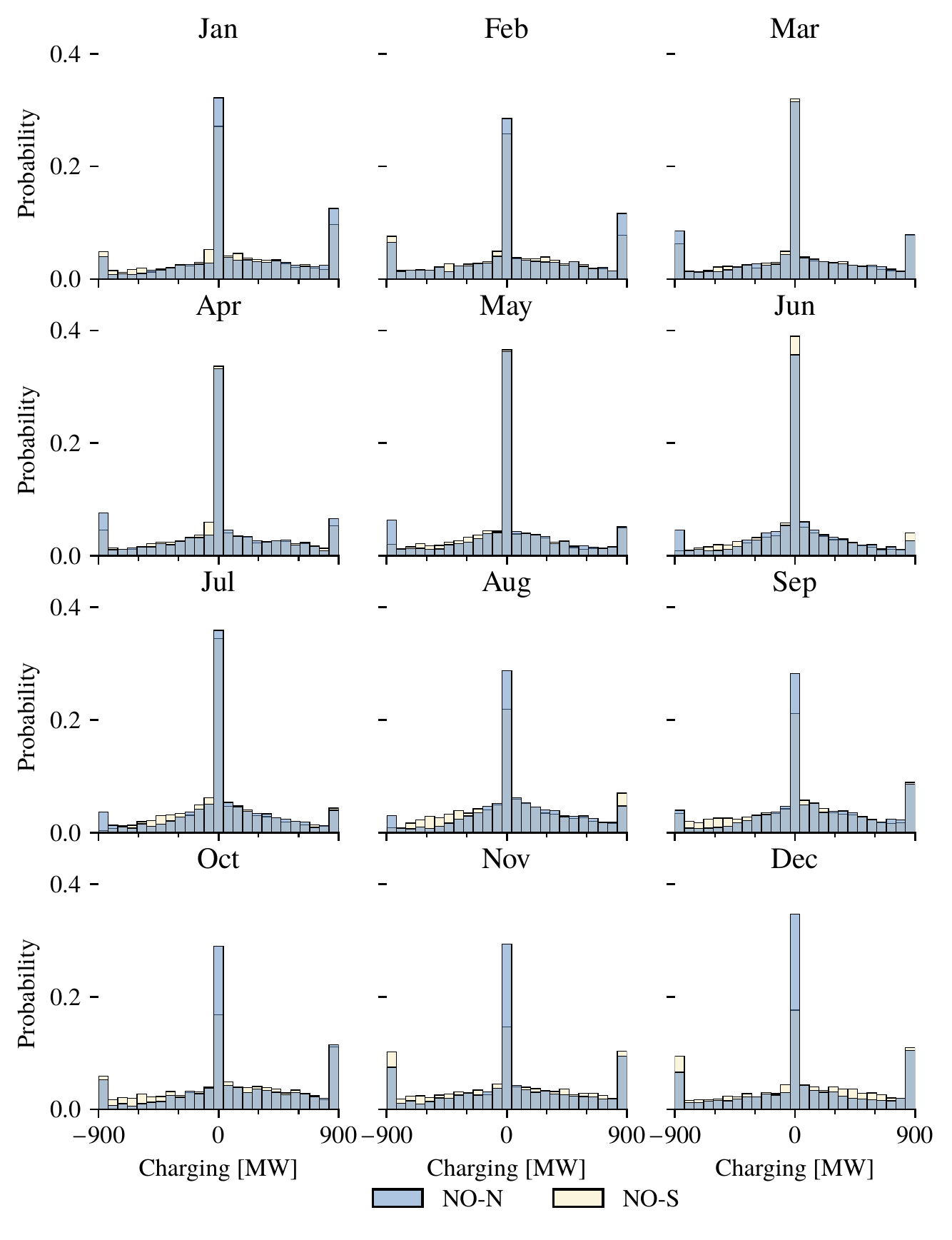}
    \caption{Distribution of storage usage throughout the different months (with 100 different calendar years) in the scenario with full flexibility.}
    \label{fig:full-flex_usage-distr}
\end{figure}

Varying the (dis-)charging capacities has a limited impact on penalty values, although they decrease with higher capacities and increase with lower capacities (as \cref{fig:stor_usage-distr} would suggest).
Higher (dis-)charging capacities (especially in NO-S) allow for higher capacities in NO-S (by up to 500~MW, thus introducing those with a ratio of 1:1) as this is a limiting factor especially for charging in the winter months (\cref{fig:stor_usage-distr} and \cref{fig:full-flex_usage-distr}).
For lower capacities, installed wind capacities are moved from south to north where conditions are better and their distribution is more favourable (\cref{fig:cap_factors}).

While the other three factors were investigated for both optimal capacities and penalty values, we only compare the latter for electricity demand, thus looking at the risk of unforeseen variations.
Consistently throughout the scenarios, higher demand increases penalty values drastically in NO-S (almost doubling for \emph{full-flex} for a uniform increase of 10\%). 
As capacities in northern Norway are rather high and prone to overproduction, increasing electricity demand actually \emph{lowers} the penalty values.
Of the parameters studied, unexpected changes in demand have the largest impact on this study. 

\section{Further graphics}
\begin{figure}[htb]
    \centering
    \includegraphics{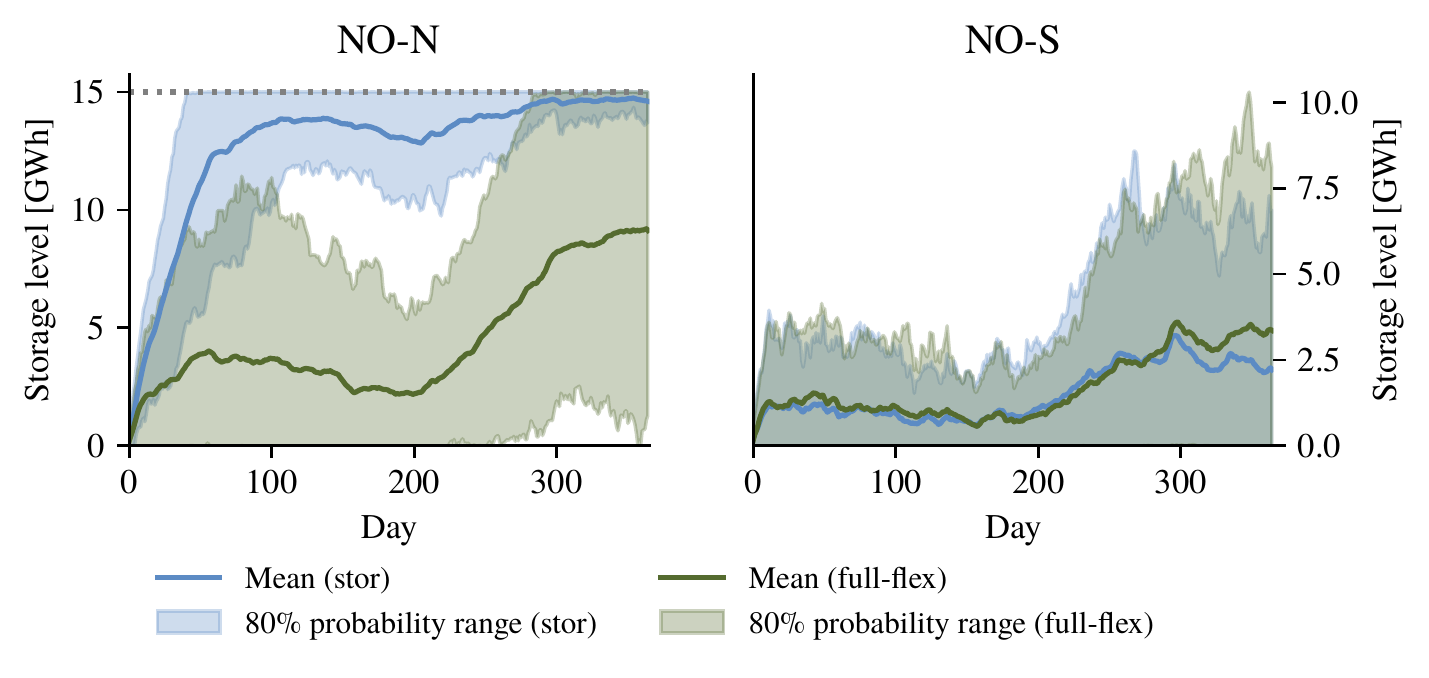}
    \caption{Comparison of mean storage level and 80\% probability ranges for the ``stor'' and ``full-flex'' scenarios. The maximal capacity for NO-N is dotted in grey for reference.}
    \label{fig:comp_level}
\end{figure}

\begin{figure}[htb]
    \centering
    \includegraphics{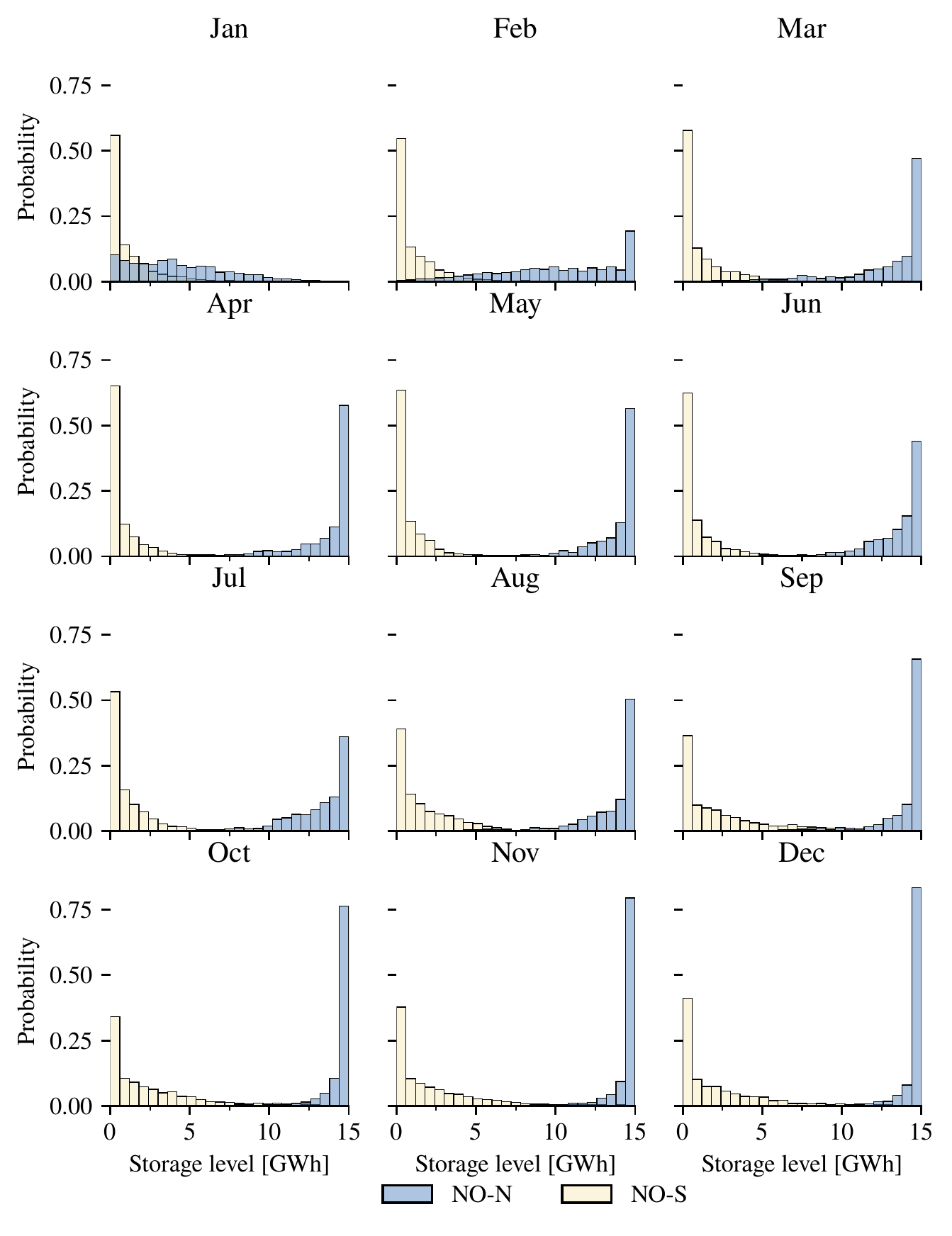}
    \caption{Distribution of storage levels throughout the different months (with 100 different calendar years) in the scenario with storage as the only flexibility option.}
    \label{fig:stor_level-distr}
\end{figure}

\begin{figure}[tbh]
    \centering
    \includegraphics{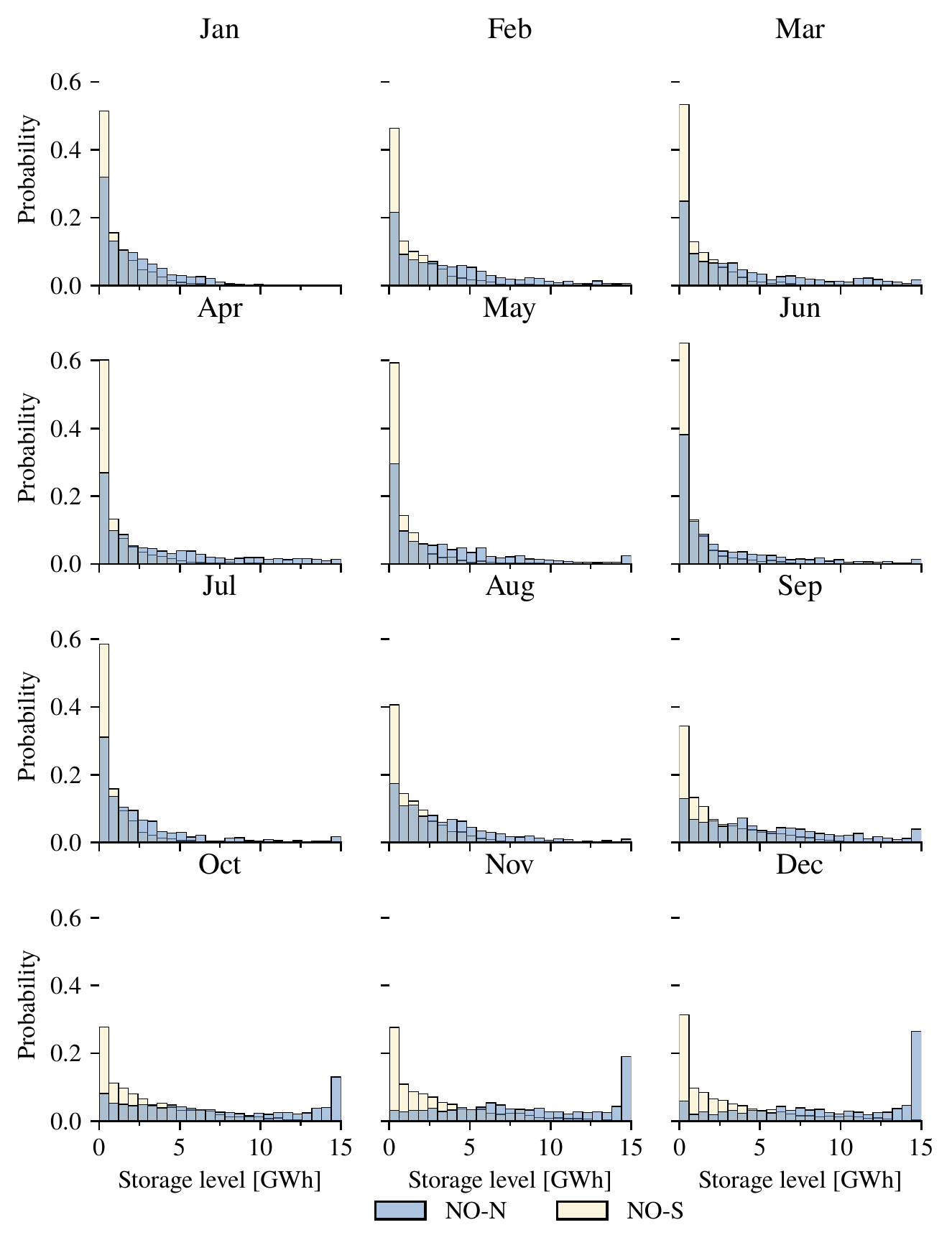}
    \caption{Distribution of storage levels throughout the different months (with 100 different calendar years) in the scenario with full flexibility.}
    \label{fig:full-flex_level-distr}
\end{figure}

\end{document}